\documentclass[runningheads]{llncs}
\usepackage[T1]{fontenc}
\usepackage[margin=1in]{geometry}
\usepackage{amsmath,amssymb,mathtools}
\usepackage{algorithm,algpseudocode,microtype}
\usepackage[colorlinks=true,allcolors=blue!55!black]{hyperref}
\usepackage{xcolor}
\usepackage{graphicx}
\microtypesetup{expansion=false}
\usepackage{tikz}
\usetikzlibrary{positioning,calc}
\usetikzlibrary{decorations.markings}
\usepackage{booktabs,tabularx,array}
\usepackage{enumitem}
\usepackage{subcaption}
\usepackage{float}
\newcommand{\HH}{\mathcal H}
\newcommand{\ArcFam}{\mathcal A}
\newcommand{\gZ}{\gamma_{\mathrm{gr}}^{Z}}
\newcommand{\ur}{\nu_{\mathrm{ur}}}

\newcommand{\EE}{\mathcal E}
\newcommand{\rgr}{\rho_{\mathrm{gr}}}
\newcommand{\tgr}{\tau_{\mathrm{gr}}}

\newcommand{\ggr}{\gamma_{\mathrm{gr}}}

\newcommand{\rh}{\rho_{\mathrm{gr}}}
\newcommand{\thh}{\tau_{\mathrm{gr}}}

\newcommand{\emptyseq}{\langle\,\rangle}
\newcommand{\cat}{\mathbin{\|}}

\begin{document}
\title{Dominating Sequence and its Variants on the Subclasses of Chordal Graphs and Co-bipartite Graphs}
\titlerunning{Dominating Sequence on Chordal Graphs and Co-bipartite Graphs}
%
\author{Suchanda Roy \and
Ramesh Hariharasubramanian\thanks{Corresponding author.}}
\authorrunning{S. Roy et al.}
%
\institute{Department of Mathematics, Indian Institute of Technology Guwahati, \\Guwahati, Assam 781039, India\\
\email{r.suchanda@iitg.ac.in}\\
\email{ramesh\_h@iitg.ac.in}}
\maketitle 
\begin{abstract}
A Grundy dominating sequence in a graph $G$ is a sequence in which each
vertex dominates at least one vertex not dominated by its predecessors,
and all vertices of $G$ are dominated at the end. The maximum length of
such a sequence is called the \emph{Grundy domination number} of $G$,
denoted by $\gamma_{gr}(G)$. Determining $\gamma_{gr}(G)$ and several
of its variants, including \emph{Total Grundy domination number}, \emph{$L$-Grundy domination number}, and
\emph{$Z$-Grundy domination number}, are computationally hard on general graphs and several restricted graph classes, including co-bipartite, split graphs. In contrast,
the Grundy domination number can be computed in polynomial time on split
graphs.

In this paper, we study the dominating sequence problem and its variants on word-representable co-bipartite graphs and circular-convex co-bipartite graphs. We obtain polynomial-time algorithms for all considered variants on both graph classes, establishing their tractability on these classes. Finally, we show that the Grundy domination number problem of a well-partitioned chordal graph is NP-complete and $Z$- Grundy domination problem as well as zero forcing set problem is NP-complete on split, bipartite and word-representable graphs.

\end{abstract}
\vspace{-1cm}

\section{Introduction}
\vspace{-.2cm}
Let \(G\) be a graph with vertex set \(V(G)\) and edge set $E(G)$. We say that a vertex $v$ \emph{dominates} a vertex \(u\) if \(u=v\) or $uv\in E(G)$.
A sequence $(v_1,v_2,\ldots,v_k)$ of vertices of a graph $G$ is called a
\emph{Grundy dominating sequence} if, for each $i\in\{1,\ldots,k\}$,
the vertex $v_i$ dominates at least one vertex that is not dominated by
any of $v_1,\ldots,v_{i-1}$, and every vertex of $G$ is dominated by the
sequence. The maximum length of a Grundy dominating sequence is called
the \emph{Grundy domination number} of $G$, denoted by
$\gamma_{gr}(G)$. The study of Grundy dominating sequences and their variants
has received considerable attention from both structural and
algorithmic perspectives, see, for example,
\cite{brevsar2023computation,brevsar2014dominating,brevsar2016dominating,nasini2020grundy,AErey}.

Several parameters for Grundy dominating sequence have been introduced by imposing additional
conditions (open or closed) on the neighborhoods of the vertices in a dominating sequence. In this
work, we consider \emph{(i) Total Grundy Dominating Sequence}
\cite{brevsar2016total,bahadir2021graphs,Dravec,lin2019zero},
\emph{(ii) $L$-Grundy Dominating Sequence}\cite{brevsarL,brevsar2017grundy,lin2019zero} and
\emph{(iii) $Z$-Grundy Dominating Sequence} \cite{brevsar2017grundy,lin2019zero}.A sequence $S=(v_1,v_2,\ldots,v_k)$ of distinct
vertices of a graph $G$ is called a
\emph{Grundy dominating sequence} if, for every
$1\leq i\leq k$,
$N[v_i]\nsubseteq \bigcup_{j=1}^{i-1}N[v_j]$,
that is, $v_i$ dominates at least one vertex
that was not dominated by its predecessors. The sequence $S$ is called a
\emph{Total Grundy dominating sequence} if
$N(v_i)\nsubseteq \bigcup_{j=1}^{i-1}N(v_j)$
for every $1\leq i\leq k$, that is, $v_i$
totally dominates a new vertex different
from itself. It is called an
\emph{$L$-Grundy dominating sequence} if
$N[v_i]\nsubseteq \bigcup_{j=1}^{i-1}N(v_j)$,
that is, $v_i$ dominates a vertex that was
not totally dominated by its predecessors.
Similarly, $S$ is called a
\emph{$Z$-Grundy dominating sequence} if
$N(v_i)\nsubseteq \bigcup_{j=1}^{i-1}N[v_j]$,
that is, $v_i$ totally dominates a vertex
that was not dominated by its predecessors.
Here, the union is empty when $i=1$. A subset $B\subseteq V(G)$ is called a
\emph{zero forcing set} if, initially
coloring all vertices of $B$ black and
the remaining vertices white, repeated
application of the following color-change
rule eventually colors every vertex black:
a black vertex with exactly one white
neighbor forces that neighbor to become
black. The minimum cardinality of a zero
forcing set is called the \emph{zero forcing
number} of $G$, denoted by $Z(G)$. Bre\v{s}ar et al.~\cite{brevsar2017grundy}
established that, for every graph \(G\) of order \(n\), the \(Z\)-Grundy domination number is equal to the number of vertices of \(G\) minus its zero forcing number. Consequently, determining the $Z$-Grundy domination
number is equivalent to determining the
zero forcing number.

Throughout, by a \emph{Grundy dominating sequence problem} (and, analogously, a
\emph{Grundy dominating sequence variant problem}) we mean the problem of finding a Grundy dominating
sequence (respectively, a Grundy dominating sequence satisfying the additional
specified constraint) of maximum length. The corresponding maximum lengths of such sequences are known as the respective \emph{Grundy domination numbers}. For a graph $G$, we use 
\(
\;\gamma_{gr}(G),\;
\gamma_{gr}^{\;t}(G),\;
\gamma_{gr}^{\;Z}(G),\;
\text{ and }
\gamma_{gr}^{\;L}(G),\;
\)
to denote \emph{Grundy domination number},
\emph{Grundy total domination number}, \emph{$Z$-Grundy domination number},
 and \emph{$L$-Grundy domination number}, respectively. The Grundy dominating sequence problem and its variants in $(i)-(iii)$ are
computationally difficult in general and remain NP-complete on several
restricted graph classes, including co-bipartite graphs\cite{scheffler2025parameterized} and chordal graphs\cite{brevsar2014dominating}\cite{scheffler2025parameterized}. In contrast, Grundy dominating sequence problem is polynomial-time solvable on split graphs and the complexity of $Z$-Grundy dominating sequence problem is still open on split and bipartite graphs. In the following, we state the result regarding Grundy domination number of split graphs.
\vspace{-.1cm}
\begin{theorem}\label{splitgrundy}\cite{brevsar2014dominating}
Let \(G=(C\cup I,E)\) be a split graph with a maximum independent set \(I\), and clique $C$. Then
\(
\gamma_{gr}(G)=
\begin{cases}
\alpha(G), & \text{if every two vertices in \(C\) }\\[.1mm]
& \text{have a common neighbor in \(I\)}\\[.2mm]
\alpha(G)+1, & \text{otherwise}.
\end{cases}
\)
\end{theorem}
\vspace{-.1cm}
Here, we connect these dominating sequence problem with the class of
\emph{word-representable graphs}. A graph is word-representable if its
adjacency relation can be characterized by alternation of letters in a
word (Definition\ref{word}). The class has been extensively studied from a structural, combinatorial and enumerative 
perspectives, leading to several characterizations and structural
results for its important subclasses; see, for example, \cite{halldorsson2011alternation,halldorsson2016semi,hefty2024word,kitaev2017comprehensive,kitaev2021word,kitaev2015words,kitaev2024semi,roy2026characterization,roy2026word,das2025word}.
It contains several familiar graph classes, including all bipartite
graphs, comparability graphs, circle graphs, 3-colorable graphs etc. Consequently, many classical graph-theoretic optimization
problems remain computationally difficult on word-representable graphs;
in particular, the dominating sequence problems and variants considered in
\emph{(i)--(ii)} are NP-complete on this class(\emph{as it is a superclass of bipartite graph and Grundy, total Grundy, $L$- Grundy dominating sequence problems are NP-complete on bipartite graphs\cite{brevsar2023computation,brevsar2016total,brevsar2017grundy}}) and for $Z$-Grundy dominating sequence problem its complexity is not known for word-representable graphs. In contrast, the algorithmic
aspects of word-representability within well-structured graph classes
have received comparatively less attention. This motivates us to
exploit the additional structure arising from word-representability to
develop efficient algorithms for various dominating sequence problems.

In particular, word-representability imposes useful structural
properties on co-bipartite graphs (Theorem \ref{coth}). Although recognizing
word-representable graphs is NP-complete \cite{kitaev2015words}, co-bipartite
graphs can be recognized in polynomial time \cite{srinivasan2026forbidden}. Building on these
structural properties,

\begin{itemize}
    \item we develop polynomial-time algorithms for the Grundy dominating sequence problem
    and the variants considered in (i)--(iii) on word-representable co-bipartite
    graphs. We further extend the structural framework to circular-convex
co-bipartite graphs (Definition \ref{codef}).  
\item On the other hand, we have shown that $Z$- Grundy domination problem as well as zero forcing set problem is NP-complete on split, bipartite and word-representable graphs and Grundy dominating sequence problem is NP-complete on well-partitioned chordal graphs.
\end{itemize}
This paper is organized as follows. Section~\ref{2} reviews the basic
definitions, results, and structural tools relevant to our work, with
a focus on word-representable graphs (Subsection \ref{2.1}) and the graph classes considered
here, namely co-bipartite (Subsection \ref{2.1.1}) and well-partitioned chordal graphs (Subsection \ref{2.2}). 
Sections~\ref{co-bip}, \ref{split} and ~\ref{wpcs} present our results
on dominating sequence and its variants for word-representable and circular-convex
co-bipartite graphs and well-partitioned
chordal graphs, respectively. Finally, Section~\ref{final} concludes the paper
with open problems and directions for future research.
\vspace{-.4cm}
\section{Preliminaries}\label{2}
\vspace{-.2cm}
We begin this section with the basic notation and definitions used
throughout the paper. All graphs considered are finite, simple, and
undirected. For a graph $G=(V,E)$, let $n=|V|$ and $m=|E|$ denote its
order and size, respectively. For vertices $x,y\in V(G)$, we write $x\sim y$ if $xy\in E$
and $x\nsim y$ otherwise. For a vertex $v\in V(G)$, its open and
closed neighborhoods are defined by
\(
N_G(v)=\{u\in V:uv\in E\}
\;\;\text{and}\;\;
N_G[v]=N(v)\cup\{v\},
\)
respectively. The degree of $v$ is denoted by $d(v)=|N_G(v)|$. For a set
$X\subseteq V$, its closed neighborhood is
\(
N_G[X]=\bigcup_{u\in X}N_G[u].
\)
The subgraph of $G$ induced by $X\subseteq V$ is denoted by $G[X]$. For a graph $G$, we denote by $G-X$ the induced subgraph
$G[V(G)\setminus X]$. For a single vertex $x\in V(G)$, we use the
shorthand $G-x$ for $G-\{x\}$. For two sets $X,Y\subseteq V(G)$, let
\(
G[X,Y]=(X\cup Y,\{xy\in E(G):x\in X,\ y\in Y\}).
\)
We say that $X$ is \emph{complete to} $Y$ if $X\cap Y=\emptyset$ and
every vertex of $X$ is adjacent to every vertex of $Y$.
For $v\in V$, we use $G-v$ to denote the graph obtained from $G$ by
deleting $v$ and all edges incident with $v$. For two sets
$S_1,S_2$, we write $S_1\setminus S_2$ for their set difference.
A set $S\subseteq V$ is called an \emph{independent set} if no two
distinct vertices of $S$ are adjacent, and a \emph{clique} if every two distinct vertices of $S$ are adjacent.

\subsection{Word-Representable Graphs}\label{2.1}
\vspace{-.1cm}
We begin this section with the definition of word-representable graphs.
\vspace{-.1cm}
\begin{definition}[\cite{kitaev2015words}]\label{word}  A simple graph $G = (V, E)$ is \textit{word-representable} if there exists a word $w$ over the alphabet $V$ such that letters $x$ and $y$, $\{x,y\}\in V$ alternate in $w$ if and only if $x\sim y$ for each $x\not=y$. If a word $w$ \textit{represents} $G$, then $w$ contains each letter of $V(G)$ at least once.
\end{definition} 
\vspace{-.1cm}
From the definition, it directly follows that, the word-representability is closed under taking induced subgraphs, i.e, it is a hereditary property.

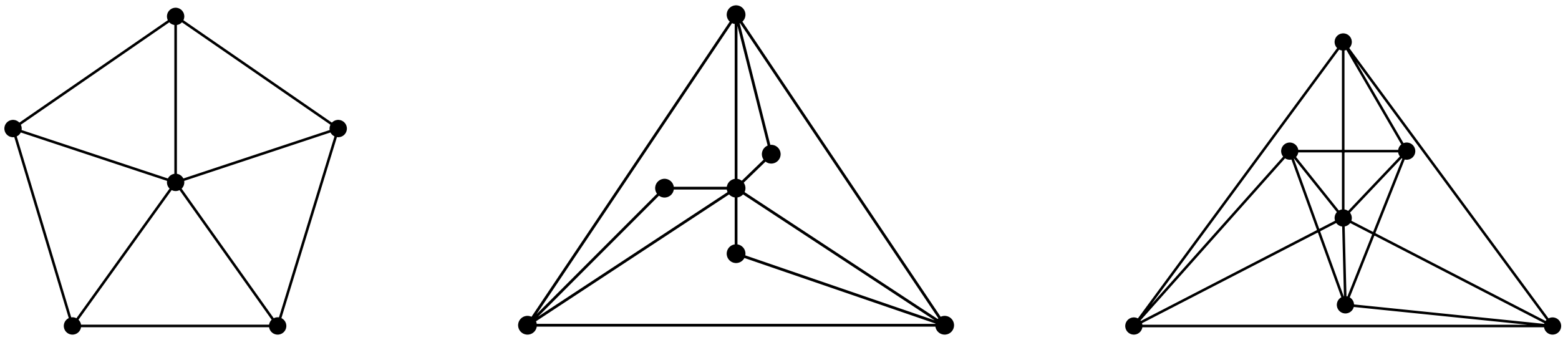
\begin{figure}[htbp]
\centering
\begin{subfigure}[t]{0.22\linewidth}
\centering
\resizebox{\linewidth}{!}{%
\begin{tikzpicture}[line cap=round,line join=round]
[xshift=4.00cm]
\coordinate (T) at (1.45,2.76);
\coordinate (L) at (0,1.76);
\coordinate (R) at (2.90,1.76);
\coordinate (BL) at (0.53,0);
\coordinate (BR) at (2.36,0);
\coordinate (C) at (1.45,1.28);
\draw[line width=0.65pt]
 (T)--(L) (T)--(R)
 (L)--(BL) (BL)--(BR) (BR)--(R)
 (T)--(C) (L)--(C) (R)--(C)
 (BL)--(C) (BR)--(C);
\fill (T) circle (2.2pt) (L) circle (2.2pt)
 (R) circle (2.2pt) (BL) circle (2.2pt)
 (BR) circle (2.2pt) (C) circle (2.2pt);

\end{tikzpicture}%
}
\caption{Smallest example of a non-word-representable graph, which is planar \cite{kitaev2015words}.}
\label{w5}
\end{subfigure}
\hfill
\begin{subfigure}[t]{0.28\linewidth}
\centering
\resizebox{\linewidth}{!}{%
\begin{tikzpicture}[line cap=round,line join=round]
[xshift=7.14cm]
\coordinate (T) at (1.72,2.56);
\coordinate (L) at (0,0);
\coordinate (R) at (3.44,0);
\coordinate (A) at (1.13,1.13);
\coordinate (B) at (1.72,1.13);
\coordinate (C) at (2.01,1.41);
\coordinate (D) at (1.72,0.59);
\draw[line width=0.65pt]
 (T)--(L) (T)--(R) (T)--(B) (T)--(C)
 (L)--(R) (L)--(A) (L)--(B)
 (R)--(B) (R)--(D)
 (A)--(B) (B)--(C) (B)--(D);
\fill (T) circle (2.2pt) (L) circle (2.2pt)
 (R) circle (2.2pt) (A) circle (2.2pt)
 (B) circle (2.2pt) (C) circle (2.2pt)
 (D) circle (2.2pt);

\end{tikzpicture}%
}
\caption{Example of a non-word-representable split graph \cite{kitaev2021word}.}
\label{sp}
\end{subfigure}
\hfill
\begin{subfigure}[t]{0.28\linewidth}
\centering
\resizebox{\linewidth}{!}{%
\begin{tikzpicture}[line cap=round,line join=round]

\coordinate (T) at (1.88,2.55);
\coordinate (L) at (0,0);
\coordinate (R) at (3.76,0);
\coordinate (A) at (1.40,1.57);
\coordinate (B) at (2.45,1.57);
\coordinate (C) at (1.88,0.97);
\coordinate (D) at (1.90,0.19);
\draw[line width=0.65pt]
 (T)--(L) (T)--(R) (T)--(C) (T)--(B)
 (L)--(R) (L)--(A) (L)--(C)
 (R)--(C) (R)--(D)
 (A)--(B) (A)--(C) (A)--(D)
 (B)--(C) (B)--(D) (C)--(D);
\fill (T) circle (2.2pt) (L) circle (2.2pt)
 (R) circle (2.2pt) (A) circle (2.2pt)
 (B) circle (2.2pt) (C) circle (2.2pt)
 (D) circle (2.2pt);

\end{tikzpicture}%
}
\caption{Example of a non-word-representable co-bipartite graph \cite{kitaev2015words}.}
\label{cobip}
\end{subfigure}
\caption{Examples of non-word-representable graphs.}
\label{fig:non-word-representable-examples}
\end{figure}
\begin{example}
For instance, the word $1234$ represents the complete graph $K_4$.
More generally, any permutation of $n$ distinct symbols represents
the complete graph $K_n$. The word $162132435465$ represents the cycle
$C_6$. In particular, several well-known graph classes, including
comparability graphs (graphs that admit a transitive orientation),
circle graphs (intersection graphs of chords of a circle),
bipartite graphs, and cycles are word-representable. However, several graph classes contain both word-representable and
non-word-representable graphs (Figures \ref{w5},\ref{sp},\ref{cobip}), including split graphs (vertex set can be partitioned into a clique and an independent set),
co-bipartite graphs (vertex set can be partitioned into
two cliques), well-partitioned chordal graphs (Definition \ref{wpc}), and planar graphs.
In this paper, our focus is on word-representable co-bipartite graphs and
well-partitioned chordal graphs. Their relevant structural properties,
definitions, and characterizations are reviewed in detail in the
following subsections.
\end{example}
\vspace{-.2cm}
\subsubsection{Word-Representability of Co-bipartite Graphs}\label{2.1.1}

A graph $G=(C_1\cup C_2,E)$ is a \emph{co-bipartite graph} if its
vertex set can be partitioned into two cliques $C_1$ and $C_2$.

An important consequence of word-representability for co-bipartite graphs is the existence of a specific ordering of the clique vertices, which plays a crucial role in our algorithmic results. In particular, word-representability imposes a circular-convex neighborhood condition on the vertices of co-bipartite graphs which is formally written in the following theorem that follows from the complete characterization of word-representable co-bipartite graphs in \cite{das2025word}.
\begin{theorem}\cite{das2025word}\label{coth}
Let $G=(C_1\cup C_2,E)$ be a word-representable co-bipartite graph, then
the vertices of each clique can be arranged in a circular order such
that, for every vertex $v\in C_1$, its neighborhood in $C_2$ forms a
circular arc in the circular order of $C_2$, and, symmetrically, for every vertex
$v\in C_2$, its neighborhood in $C_1$ forms a circular arc in the circular order of $C_1$.    
\end{theorem}

In this context, we define circular-convex co-bipartite graphs.
\begin{definition}\label{codef}
 Let $G=(A\mathbin{\dot\cup}B,E)$ be co-bipartite, that is, $A$ and $B$ each induce a clique.  We call $G$ \emph{circular-convex with respect to $A$} if $A$ admits a cyclic ordering in which $N(b)\cap A$ is a circular interval for every $b\in B$.   
\end{definition}

In general, by  circular-convex co-bipartite graph we mean a co-bipartite graph which is circular-convex with respect to atleast one clique. Throughout this work, we denote the classes of word-representable co-bipartite graphs and circular-convex co-bipartite graphs by $\mathcal{WR}_{\mathrm{co\text{-}bip}}$ and $\mathcal{CC}_{\mathrm{co\text{-}bip}}$, respectively.
Hence, from Theorem \ref{coth}, we have the hierarchy:
\(
\mathcal{WR}_{\mathrm{co-bip}}
\subseteq
\mathcal{CA}_{\mathrm{co-bip}}
\subseteq
\small{Co\text{-}bipartite}\\ \small{Graphs}.
\)
 Consequently, any problem that admits
a polynomial-time algorithm on $\mathcal{CC}_{\mathrm{co-bip}}$ is
immediately polynomial-time solvable on
$\mathcal{WR}_{\mathrm{co-bip}}$.

\subsection{Well-Partitioned Chordal Graphs}\label{2.2}

\begin{definition}\label{wpc}\cite{ahn2022well}
A connected graph $G$ is called a \emph{well-partitioned chordal graph}
if there exists a partition $\mathcal{P}$ of $V(G)$ and a tree $T$ with
$V(T)=\mathcal{P}$ such that the following conditions hold:
\begin{enumerate}
\item[(i)] Every part $X\in\mathcal{P}$ induces a clique in $G$.

\item[(ii)] For every edge $XY\in E(T)$, there exist subsets
$X'\subseteq X$ and $Y'\subseteq Y$ such that
\(
    E(G[X,Y])\\=\{xy\mid x\in X',\,y\in Y'\},
\)
where $X'$ is complete to $Y'$. We call $X'$ the
\emph{pocket} of $X$ with respect to $Y$, denoted by
$\operatorname{p}(X,Y)$, and call $Y'$ as $\operatorname{p}(Y,X)$
analogously.

\item[(iii)] For every pair of distinct bags $X,Y\in V(T)$ with
$XY\notin E(T)$,
\(
    E(G[X,Y])=\emptyset.
\)

\end{enumerate}
The tree $T$ is called a \emph{partition tree} of $G$, and the
elements of $\mathcal{P}$ are called \emph{bags}. A graph is a
well-partitioned chordal graph if each of its connected components is a
well-partitioned chordal graph.
\end{definition}
In the following, we review some results regarding legal covering sequence of hypergraphs.
\subsection{Legal Covering Sequence of hypergraph}

Let $\HH=(X,\EE)$ be a hypergraph. We distinguish indexed
hyperedges even when two of them have the same vertex set.  A sequence
$(E_1,\ldots,E_k)$ of distinct indexed hyperedges is a \emph{legal
covering sequence} if
\begin{equation}\label{eq:legalcover}
 E_i\setminus\bigcup_{j<i}E_j\ne\varnothing\qquad(i\in[k]).
\end{equation}
The maximum length is the \emph{Grundy covering number} $\rgr(\HH)$,
introduced and studied by Bre\v{s}ar et al.~\cite{brevsar2014dominating}.
An $x_i\in E_i\setminus\bigcup_{j<i}E_j$ is a \emph{fresh witness}.

A \emph{legal transversal sequence} is an ordered sequence of
points $(x_1,\ldots,x_k)$ such that, for each $i$, some hyperedge
$E_i$ contains $x_i$ but none of $x_1,\ldots,x_{i-1}$.  Its maximum
length is $\tgr(\HH)$.  Empty hyperedges never contribute to either
parameter.

\begin{lemma}\cite{brevsar2016total}\label{lem:duality}
For every hypergraph $\HH$,
\begin{equation}\label{eq:duality}
 \rgr(\HH)=\tgr(\HH).
\end{equation}
\end{lemma}

A hypergraph is an \emph{interval hypergraph} if its points have a
linear order in which each nonempty hyperedge is consecutive.  In a
\emph{circular-arc hypergraph}, an analogous ordering is cyclic.

A \emph{max point-tolerance} (MPT) representation assigns to every
vertex $v$ an interval $I_v$ and a distinguished point $p_v\in I_v$.
Two distinct vertices $u,v$ are adjacent if and only if both
$p_u,p_v$ lie in $I_u\cap I_v$.  Catanzaro et al.
\cite{catanzaro2017max} proved that maximum-weight independent set is
polynomial-time solvable on MPT graphs when such a representation is
provided. We will be using this algorithm as a subroutine in our results.


In the following, we review some results regarding uniquely restricted matching on bipartite graphs.
\subsection{Uniquely Restricted Matching Characterization and Complexity}

A matching $M$ in a graph $H$ is \emph{uniquely restricted} if it is the
unique perfect matching of the subgraph induced by its saturated
vertices.  Write $\ur(H)$ for its maximum cardinality.
Golumbic, Hirst and Lewenstein~\cite{golumbic2001uniquely}
proved that maximizing this cardinality is NP-complete on bipartite
graphs and described the alternating-cycle characterization.

\begin{lemma}
\cite{golumbic2001uniquely}\label{lem:UR}
Let $H=(U,W;E_H)$ be bipartite and let
$M=\{u_1w_1,\ldots,u_rw_r\}$ be a matching.
Then $M$ is uniquely restricted if and only if its edges can be
ordered such that
\begin{equation}\label{eq:triangular}
 u_iw_j\notin E_H\qquad\text{whenever }i<j.
\end{equation}
\end{lemma}

\section{Dominating Sequence and its Variants on Co-Bipartite Graphs}\label{co-bip}
\vspace{-.25cm}
In this section, we first study the Maximum Grundy Dominating Sequence problem and give a polynomial-time algorithm for word-representable co-bipartite as well as circular-convex co-bipartite graphs. We then use it to develop polynomial-time algorithms for the other variants considered in this work. Before this, we strudy legal covering sequences in interval and circular-Arc hypergraphs which we will be using in our algorithms.
\subsection{Legal Covering Sequences in Interval and Circular-Arc Hypergraphs}
\label{sec:covering}



Let $\HH=(X,\EE)$ be an interval hypergraph with a fixed linear
ordering of $X$. We define an auxiliary graph $Q_{\HH}$ whose vertices
are the hypergraph incidences:
\(
 V(Q_{\HH})=\{(x,E): E\in\EE,\ x\in E\}.
\)
Two distinct incidences are adjacent precisely when
 \begin{equation}\label{eq:mptadj}
(x,E)(y,F)\in E(Q_{\HH})
\quad\Longleftrightarrow\quad
x,y\in E\cap F.
\end{equation}

We regard $(x,E)$ as the pointed interval $(E,x)$ and since every $E$ is an interval in the same linear order,
\eqref{eq:mptadj} gives an explicit MPT representation of
$Q_{\HH}$.

\begin{lemma}\label{lem:interval-id}
For every interval hypergraph $\HH$,
 $\rgr(\HH)=\alpha(Q_{\HH}),$
where $\alpha$ is the maximum independent-set cardinality.
\end{lemma}
\begin{proof}
We prove the two inequalities separately.

\emph{From a legal covering sequence to an independent set:}
Let $(E_1,\ldots,E_k)$ be legal and choose a fresh witness $x_i$ for
each selection.  If $i<j$, then $x_j\notin E_i$, since $x_j$ was not
covered when $E_j$ was selected.  Thus the incidence vertices
$(x_i,E_i)$ and $(x_j,E_j)$ are nonadjacent by \eqref{eq:mptadj}.
The $k$ incidences form an independent set, and hence
$\rgr(\HH)\le\alpha(Q_{\HH})$.

\emph{From an independent set to a legal covering sequence:}
Now suppose $\{(x_i,E_i):i\in[k]\}$ is an independent set.
Its pointed positions are pairwise distinct, since incidences sharing
one point are adjacent.  Similarly, its indexed hyperedges are
distinct. We construct a digraph $D$ on $[k]$ by placing an arc
\(
 i\longrightarrow j \quad\text{if }i\ne j\text{ and }x_j\in E_i.
\)
Independence prevents a directed two-cycle: if $i\to j$ and $j\to i$,
then $x_i,x_j$ lie in both $E_i$ and $E_j$, which would make the two
incidences adjacent.

We claim that $D$ has no longer directed cycle either.  Suppose, to
the contrary, that it has a shortest directed cycle. Then choose a vertex
$i$ on this cycle whose pointed position $x_i$ is leftmost, and let
$k$ and $j$ be its predecessor and successor.  Thus
$x_i\in E_k$ and $x_j\in E_i$.  Both $x_k$ and $x_j$ lie strictly to
the right of $x_i$.
If $x_k<x_j$, the interval $E_i$ contains $x_k$, so $i\to k$,
contradicting the absence of two-cycles.  If $x_j<x_k$, the interval
$E_k$, containing $x_i$ and $x_k$, also contains $x_j$.  The arc
$k\to j$ then shortcuts the chosen cycle, contradicting its
minimality.  Hence $D$ is acyclic.

Next, we take a reverse topological ordering $i_1,\ldots,i_k$ of $D$.
For any $s<t$, the earlier hyperedge $E_{i_s}$ cannot contain the
later pointed position $x_{i_t}$: that would create an arc
$i_s\to i_t$, which is impossible in a reverse topological order.
Consequently, $E_{i_t}$ newly covers $x_{i_t}$ when selected.
The sequence $(E_{i_1},\ldots,E_{i_k})$ is legal and has length $k$.
This proves the reverse inequality and establishes $\rgr(\HH)=\alpha(Q_{\HH})$.
\end{proof}

\begin{theorem}\label{thm:interval}
Given a linear interval representation of $\HH$, one can compute
$\rgr(\HH)$ and construct an optimal legal covering sequence in
polynomial time.
\end{theorem}
\begin{proof}
The graph $Q_{\HH}$ has at most $|X|\,|\EE|$ vertices and an explicit
pointed-interval representation.  By Lemma~\ref{lem:interval-id}, its
independence number equals $\rgr(\HH)$.  Now, we apply the polynomial-time
maximum-weight independent-set algorithm for represented MPT graphs
\cite{catanzaro2017max} with unit weights.  Then the reverse-topological
construction in the preceding proof recovers an optimal legal sequence.
\end{proof}

For a hypergraph $\HH=(X,\EE)$ and $x\in X$, let
 \(\EE_x=\{E\in\EE:x\notin E\},\qquad
 \HH_x=(X\setminus\{x\},\EE_x).\)
Note that this construction deletes \emph{all} hyperedges that
contain $x$.

\begin{lemma}\label{lem:last}
If some indexed hyperedge of $\HH$ is nonempty, then
 $\rgr(\HH)=1+
 \max_{\substack{x\in Xx\in E\\\text{ for some }E\in\EE}}
 \rgr(\HH_x).$
If all hyperedges are empty, then $\rgr(\HH)=0$.
\end{lemma}
\begin{proof}
First we take a maximum legal covering sequence $(E_1,\ldots,E_k)$.
Because its last selection is legal, it has a point
$x\in E_k\setminus\bigcup_{i<k}E_i$.
Every earlier hyperedge avoids $x$, so the first $k-1$ hyperedges
form a legal sequence of $\HH_x$.  Hence
$k\le1+\rgr(\HH_x)$.

For the reverse inequality, we fix any point $x$ lying in a nonempty
hyperedge $F$. Then take a maximum legal sequence of $\HH_x$ and append
$F$.  None of the preceding hyperedges contains $x$, whereas $F$
contains $x$.  Thus $F$ is legal as the last selection, and the new
sequence has length $\rgr(\HH_x)+1$.  Now, maximizing over $x$ gives the required result. The empty-hyperedge case is immediate.
\end{proof}

The above identity holds for arbitrary hypergraphs. The useful
feature of a \emph{circular-arc} hypergraph is that every hyperedge
avoiding $x$ becomes an ordinary interval after the cyclic order is
cut at $x$.

\begin{algorithm}[t]
\caption{Maximum legal covering sequence of circular arcs}
\label{alg:circular}
\begin{algorithmic}[1]
\Require Circular-arc hypergraph $\HH=(X,\ArcFam)$ and a cyclic order of $X$
\If{every indexed arc is empty}
  \State \Return $(0,\emptyseq)$
\EndIf
\State $\mathrm{best}\gets0$; $S_{\mathrm{best}}\gets\emptyseq$
\For{each $x\in X$ contained in some indexed arc}
 \State $\ArcFam_x\gets\{F\in\ArcFam:x\notin F\}$
 \State Cut the cyclic ordering at $x$ to represent the arcs in $\ArcFam_x$ as intervals
 \State Compute an optimal legal sequence $S_x$ of $\HH_x$ using Theorem~\ref{thm:interval}
 \State Choose an indexed arc $F_x\in\ArcFam$ containing $x$
 \If{$|S_x|+1>\mathrm{best}$}
  \State $\mathrm{best}\gets|S_x|+1$
  \State $S_{\mathrm{best}}\gets S_x\cat\langle F_x\rangle$
 \EndIf
\EndFor
\State \Return $(\mathrm{best},S_{\mathrm{best}})$
\end{algorithmic}
\end{algorithm}

\begin{theorem}\label{thm:circular}
The maximum legal covering-sequence problem is polynomial-time
solvable for circular-arc hypergraphs with a supplied cyclic ordering.
\end{theorem}
\begin{proof}
Algorithm~\ref{alg:circular} enumerates a possible last fresh point
$x$.  All hyperedges remaining in $\HH_x$ avoid $x$ and hence are
intervals when the circle is cut at $x$.  Theorem~\ref{thm:interval}
computes their optimum legal sequence.  Appending an arc containing
$x$ gives a valid circular-arc sequence.  By Lemma~\ref{lem:last},
one of the choices of $x$ achieves the optimum.

If $m=|X|$ and $p=|\ArcFam|$, there are at most $m$ interval instances,
each with at most $m$ points and $p$ indexed arcs.  If
$T_{\mathrm{MPT}}(N)$ is a polynomial bound for the represented-MPT
independent-set algorithm on $N$ vertices, a direct implementation
runs in
 $O\!\left(m\bigl((mp)^2+T_{\mathrm{MPT}}(mp)\bigr)\right).$
In particular, the procedure is polynomial. 
\end{proof}

\subsection{Ordinary and $Z$-Grundy Domination on Circular-convex Co-bipartite Graphs}
We first introduce the notion of \emph{cross-incidence hypergraph} hypergraph of a co-bipartite graph. For a co-bipartite graph $G=(A\mathbin{\dot\cup}B,E)$ 
the cross-incidence hypergraph of $G$ with respect
to $A$ is the hypergraph
$\mathcal H_A=(A,\mathcal E_A),$
where the ground set is $A$ and
$\mathcal E_A=
\{E_b=N_G(b)\cap A:
b\in B,\ N_G(b)\cap A\neq\varnothing\}.$ 
Each hyperedge $E_b$ is indexed by its corresponding
vertex $b\in B$. Thus, even if two distinct vertices
$b,b'\in B$ have identical neighborhoods in $A$,
their corresponding hyperedges are treated as distinct
indexed hyperedges.
Similarly, the cross-incidence hypergraph of $G$
with respect to $B$ is
$\mathcal H_B=(B,\mathcal E_B),$
where
$\mathcal E_B=
\{E_a=N_G(a)\cap B:
a\in A,\ N_G(a)\cap B\neq\varnothing\}.$
\begin{lemma}
Let $\mathcal H_A$ and $\mathcal H_B$ denote the
cross-incidence hypergraphs of $G$ with respect
to $A$ and $B$, respectively. Then $\rh(\mathcal H_A)
=\rh(\mathcal H_B).$
\end{lemma}
\begin{proof}
 Since these hypergraphs are incidence duals,
the covering--transversal duality implies
$\rh(\mathcal H_A)
=\thh(\mathcal H_A)
=\rh(\mathcal H_B)
=\thh(\mathcal H_B).$
Thus the maximum length of a legal
cross-neighborhood sequence is the same
regardless of whether the selected vertices
belong to $A$ or $B$.

Consequently, it is sufficient to compute
this parameter using $\mathcal H_A$.
If $G$ is circular-convex with respect to $A$,
then $\mathcal H_A$ is a circular-arc hypergraph,
and its legal covering number can be computed
in polynomial time.
\end{proof}
\begin{theorem}\label{thm:grz}
For every connected co-bipartite graph $G=(A\mathbin{\dot\cup}B,E)$ of order at least two, with $A,B\ne\varnothing$,
$\ggr(G)=\gZ(G)=t+\varepsilon_A+\varepsilon_B,$ where $t=\rh(\mathcal H_A)=\thh(\mathcal H_A)$ and $\varepsilon_P=\begin{cases}1,&\exists p\in P:N(p)\cap Q=\varnothing,\\0,&\text{otherwise},\end{cases}$
where $\{P,Q\}=\{A,B\}$.
In particular, both parameters are polynomial-time computable on circular-convex co-bipartite graphs.
\end{theorem}
\begin{proof}
Consider an ordinary Grundy sequence starting in $A$.  Before it first selects a vertex of $B$, only the first selected $A$-vertex may be legal without a fresh cross-neighbor (the neighbors in the other partition), as after that first selection, all of $A$ is closed-dominated.  Hence its $A$-prefix has at most $t+\varepsilon_A$ vertices.  The first $B$-selection dominates all remaining vertices of $B$, and then the sequence ends.  If $\varepsilon_B=0$, a legal final $B$-selection has an undominated witness $b\in B$ with a cross-neighbor in $A$.  Some unused $A$-vertex could instead extend the cross-neighborhood sequence to dominate $b$, showing that the extra $B$-selection cannot increase the total bound beyond $t+\varepsilon_A$.  If $\varepsilon_B=1$, at most one extra $B$-selection is possible.  The symmetric argument applies to sequences starting in $B$.  Consequently $\ggr(G)\le t+\varepsilon_A+\varepsilon_B$.

For the reverse inequality, we take a legal cross-neighborhood sequence in $A$ of length $t$; each selected vertex has a fresh neighbor in $B$, so it is $Z$-legal.  If $\varepsilon_A=1$, we prepend a cross-isolated vertex of $A$, which has a neighbor within $A$ by connectedness and is therefore $Z$-legal as first selection.  If $\varepsilon_B=1$, we choose a cross-isolated $b_0\in B$ and append a different $b\in B$ as it newly open-dominates $b_0$, which none of the preceding $A$-vertices closed-dominates.  This yields a $Z$-Grundy sequence of the asserted length.  Since $\gZ(G)\le\ggr(G)$, equality follows.

One-sided circular convexity makes $\mathcal H$ a circular-convex hypergraph, so Theorem~\ref{thm:circular} computes $t$ in polynomial time.
\end{proof}

\subsection{Total and $L$-Grundy domination on Circular-convex Co-bipartite Graphs}
For $\sigma\in\{t,L\}$, let $F_t(v)=N(v)$ and $F_L(v)=N[v]$. In both variants, before selecting $v_i$,
the vertices already dominated are those in
$\bigcup_{j=1}^{i-1}N(v_j).$

\begin{lemma}\label{lem:suffix}
Suppose a total Grundy or $L$-Grundy sequence starts in $A$, and $b_1$ is the first selected vertex in $B$.  Then at most two additional vertices can follow $b_1$.  The sequence therefore has the form $A^*\cdot T$, where $T$ is empty or starts in $B$, with $|T|\le3$.
\end{lemma}
\begin{proof}
Let $a_1$ be the first selected vertex of $A$.  The union $N(a_1)\cup N(b_1)$ contains every vertex except possibly $a_1,b_1$.  Every later selection must have a fresh footprint outside the union of earlier open neighborhoods, and there are at most two candidates.  Thus at most two later selections are possible.
\end{proof}
\begin{lemma}\label{lem:prefix}
If a total Grundy or $L$-Grundy sequence starts with two distinct vertices $a_1,a_2\in A$, then every further selection in $A$ before the first selection in $B$ must newly open-dominate a vertex of $B$.
\end{lemma}
\begin{proof}
Since $A$ is a clique, $N(a_1)\cup N(a_2)$ covers all of $A$.  Subsequent $A$-selections can therefore have fresh footprints only in $B$.
\end{proof}

\begin{lemma}\label{lem:restricted}
Let $G=(A\mathbin{\dot\cup}B,E)$ be a co-bipartite
graph that is circular-convex with respect to $A$. For any $P\subseteq A$ and $Q\subseteq B$, we can
compute, in polynomial time, the maximum length of
a sequence of vertices from $P$ such that each
selected vertex newly dominates at least one
vertex of $Q$. The same result holds when $P\subseteq B$ and
$Q\subseteq A$.
\end{lemma}

\begin{proof}
Without loss of generality, let $P\subseteq A$ and $Q\subseteq B$.

We first consider the restricted cross-incidence hypergraph
\(\mathcal H_{P,Q}=(P,\mathcal E_{P,Q})\),
where
\(\mathcal E_{P,Q}
=\{N_G(q)\cap P:q\in Q,\
N_G(q)\cap P\neq\varnothing\}\).
A sequence of vertices from $P$, in which each
selected vertex newly dominates a vertex of $Q$,
corresponds exactly to a legal transversal
sequence of $\mathcal H_{P,Q}$.
Therefore, the maximum length of such a sequence
is $\thh(\mathcal H_{P,Q})$.
By the covering--transversal duality,
\(\thh(\mathcal H_{P,Q})
=\rh(\mathcal H_{P,Q})\).
Since $G$ is circular-convex with respect to $A$,
each hyperedge $N_G(q)\cap P$ is a circular
interval in the cyclic ordering restricted to $P$.
Hence, by Theorem~\ref{thm:circular},
$\rh(\mathcal H_{P,Q})$ can be computed in
polynomial time.

The same argument holds when $P\subseteq B$ and
$Q\subseteq A$. Hence, in both cases, the required maximum
length can be computed in polynomial time.
\end{proof}

\begin{algorithm}[ht]
\caption{Total or $L$-Grundy domination on circular-convex co-bipartite graphs}\label{alg:totalL}
\begin{algorithmic}[1]
\Require A circular-convex co-bipartite graph $G=(A\mathbin{\dot\cup}B,E)$ and $\sigma\in\{t,L\}$.
\State $\texttt{best}\gets0$.
\For{each ordered choice $(P,Q)\in\{(A,B),(B,A)\}$}
  \For{each $p_1\in P$ and each ordered $T$ of at most three other vertices}
    \If{$T$ is empty or its first member lies in $Q$}
      \State Check $(p_1)\cdot T$ directly and update $\texttt{best}$ if legal.
    \EndIf
  \EndFor
  \For{each ordered pair of distinct $p_1,p_2\in P$}
    \For{each ordered suffix $T=(t_1,\ldots,t_s)$, $s\le3$, disjoint from $\{p_1,p_2\}$, empty or starting in $Q$}
      \For{each witness tuple $(w_1,\ldots,w_s)\in Q^s$}
        \If{$w_i\in F_\sigma(t_i)\setminus\bigl(N(p_1)\cup N(p_2)\cup\bigcup_{j<i}N(t_j)\bigr)$ for all $i$}
          \State $W\gets\{w_1,\ldots,w_s\}$.
          \State $P'\gets P\setminus\bigl(\{p_1,p_2\}\cup(T\cap P)\cup N(W)\bigr)$.
          \State $Q'\gets Q\setminus\bigl(N(p_1)\cup N(p_2)\bigr)$.
          \State $q\gets$ maximum legal cross-neighborhood length from $P'$ into $Q'$.
          \State $\texttt{best}\gets\max\{\texttt{best},2+q+s\}$.
        \EndIf
      \EndFor
    \EndFor
  \EndFor
\EndFor
\State \Return $\texttt{best}$.
\end{algorithmic}
\end{algorithm}

\begin{theorem}\label{thm:totalL}
Both Maximum Total Grundy Dominating Sequence and Maximum $L$-Grundy Dominating Sequence are polynomial-time solvable on circular-convex co-bipartite graphs.
\end{theorem}
\begin{proof}
We consider an optimal $\sigma$-Grundy sequence and choose $(P,Q)$ so that its first selected vertex lies in $P$.  By Lemma~\ref{lem:suffix}, the suffix after the first switch has at most three vertices.  If the initial $P$-prefix has length one, Algorithm~\ref{alg:totalL} explicitly tests the entire sequence.

Otherwise, the optimum can be written as $(p_1,p_2)\cdot R\cdot T$, where $R$ is a sequence of further vertices of $P$ and $T$ is empty or starts in $Q$. For each $t_i\in T$, a genuine fresh $\sigma$-footprint witness $w_i$ can be selected.  Since the first two selected $P$-vertices have open-dominated all of $P$, every $w_i$ lies in $Q$.  The guessed witnesses also avoid the open neighborhoods of $p_1,p_2$ and the preceding suffix vertices, and no member of $R$ may dominate any witness.  Moreover, $R$ is disjoint from $T$.  Thus $R$ lies in the algorithm's $P'$ and forms a legal cross-neighborhood sequence into $Q'$, by Lemma~\ref{lem:prefix}.  Consequently, $|R|\le q$, and the algorithm considers a candidate at least as long as the optimum.

Conversely, for every accepted enumeration, which begins with the two clique vertices $p_1,p_2$; both are legal for total and $L$ since $P$ is a clique. Then, our algorithm appends a maximum legal cross-neighborhood sequence of length $q$ using $P'$.  Its selections remain legal because they introduce fresh vertices of $Q'$.  The witness restrictions guarantee that each $w_i$ remains undominated before the corresponding $t_i$, even after the cross-neighborhood prefix. Hence, the enumerated suffix $T$ is legal as well. Thus, every candidate corresponds to a feasible sequence.  The algorithm returns the optimum.

There are polynomially many choices of the two initial vertices, the at-most-three suffix vertices, and the at-most-three footprint witnesses.  Each calls the polynomial subroutine in Lemma~\ref{lem:restricted}. Hence, the overall algorithm is polynomial.
\end{proof}

\begin{corollary}
For circular-convex co-bipartite graph, ordinary Grundy domination, $Z$-Grundy domination, total Grundy domination, and $L$-Grundy domination are all polynomial-time solvable.
\end{corollary}
\section{NP-completeness of $Z$-Grundy Domination on Split and Bipartite Graphs}
\label{split}

We prove both hardness results through the same source problem. Let $H=(A \cup B, E(H))$ be a bipartite graph. From $H$, we construct two new graphs $G_H^{\mathrm{sp}}$ and $G_H^{\mathrm{bi}}$. For every $a\in A$, we introduce two distinct vertices $a^0,a^1$ and put
 $A^*=\{a^0,a^1:a\in A\}.$
The \emph{common core} has vertex set $A^*\cup B$ and precisely the
edges
\begin{equation}\label{eq:zcore}
 a^\ell b\in E\quad\Longleftrightarrow\quad ab\in E_H,
 \qquad \ell\in\{0,1\}.
\end{equation}
Thus $a^0$ and $a^1$ have the same neighborhood in $B$.
We now describe two ways to complete this core.

\smallskip\noindent
\emph{Split construction $G_H^{\mathrm{sp}}$.}
Make $B\cup\{c_0,c_1\}$ a clique and $A^*$ an independent set.  The
vertices $c_0,c_1$ have no neighbors in $A^*$; the cross-edges are
exactly those in~\eqref{eq:zcore}.  The resulting graph is a connected
split graph.

\smallskip\noindent
\emph{Bipartite construction $G_H^{\mathrm{bi}}$.}
Keep the core bipartite, add a vertex $h$ to the $A^*$ side and two
vertices $p_0,p_1$ to the $B$ side, and add precisely $\{hb:b\in B\}\cup\{hp_0,hp_1\}.$
Thus
$ V(G_H^{\mathrm{bi}})=
 (A^*\cup\{h\})\mathbin{\dot\cup}(B\cup\{p_0,p_1\}).$
Every vertex of $B$ is adjacent to $h$, both $p_0,p_1$ are leaves at
$h$, and no other edges are present.  This graph is also connected.
Both constructions have size polynomial in $|V(H)|+|E_H|$.


\begin{lemma}\label{lem:zcorebounds}
Let $T$ be any $Z$-Grundy sequence in $G_H^{\mathrm{sp}}$ or
$G_H^{\mathrm{bi}}$, and $r=\nu_{\mathrm{ur}}(H)$.
Then:
\begin{enumerate}
\item at most $r$ selected vertices belong to $A^*$;
\item at most $r$ selected vertices belong to $B$ and have a fresh
      footprint in $A^*$.
\end{enumerate}
\end{lemma}
\begin{proof}
Let $a_1^{\ell_1},\ldots,a_s^{\ell_s}$ be the list of the selected vertices of $A^*$, in their order of appearance.
  Each must have a fresh footprint
$b_i\in B$, since its entire open neighborhood lies in $B$.  Two
copies of the same $a\in A$ cannot both be selected as the first copy
closed-dominates all possible footprints of the second.  Also the
footprints $b_i$ are distinct. Also, freshness of each $b_i$ gives
$a_ib_j\notin E_H$ whenever $i<j$.  Hence
$\{a_1b_1,\ldots,a_sb_s\}$ is a uniquely restricted matching by
Lemma~\ref{lem:UR}, so $s\le r$.

Next, we list the selected vertices of $B$ that have a fresh footprint
in $A^*$ as $b_1',\ldots,b_t'$, again in their order of appearance.
Next, we choose a fresh footprint $a_i^{\ell_i}$ for each $b_i'$.  The
underlying vertices $a_i\in A$ are distinct; otherwise selecting a neighbor of
one copy of $a$ closed-dominates both copies.  Since an earlier
$b_i'$ cannot be adjacent to the later fresh footprint
$a_j^{\ell_j}$, we have $a_jb_i'\notin E_H$ for $i<j$.
Reversing these matching pairs gives the triangular ordering of
Lemma~\ref{lem:UR}.  Thus $t\le r$.  
\end{proof}

\begin{theorem}\label{thm:zidentities}
For every bipartite graph $H=(A\mathbin{\dot\cup}B,E_H)$ without
isolated vertices and with at least one edge, the above constructions
satisfy
\begin{align}
 \gamma_{\mathrm{gr}}^Z(G_H^{\mathrm{sp}})
    &=2\nu_{\mathrm{ur}}(H)+1,\label{eq:zspidentity}\\
 \gamma_{\mathrm{gr}}^Z(G_H^{\mathrm{bi}})
    &=2\nu_{\mathrm{ur}}(H)+2.\label{eq:zbiidentity}
\end{align}
\end{theorem}
\begin{proof}
Write $r=\nu_{\mathrm{ur}}(H)$ and choose a matching
$M=\{a_1b_1,\ldots,a_rb_r\}$ ordered as in
Lemma~\ref{lem:UR}.

\smallskip\noindent
\emph{Lower bound for the split graph.}
Consider $S_{\mathrm{sp}}=(a_1^0,\ldots,a_r^0,c_0,b_r,\ldots,b_1).$
For each $i$, the vertex $a_i^0$ has fresh footprint $b_i$, by
\eqref{eq:triangular}.  The vertex $c_0$ has fresh footprint $c_1$,
since none of the preceding copied vertices is adjacent to $c_1$.
Finally, in the reverse suffix, $b_i$ has fresh footprint $a_i^1$:
this vertex is not adjacent to $c_0$, and no previously selected
$b_j$ with $j>i$ is adjacent to it.  Thus
$S_{\mathrm{sp}}$ is $Z$-Grundy and has length $2r+1$.

\smallskip\noindent
\emph{Lower bound for the bipartite graph.}
Here use $S_{\mathrm{bi}}=(p_0,a_1^0,\ldots,a_r^0,h,b_r,\ldots,b_1).$
The first vertex $p_0$ footprints $h$; each $a_i^0$ footprints $b_i$;
$h$ footprints $p_1$; and each $b_i$
footprints $a_i^1$, as in the split case.  In particular, the new
edges at $h$ do not reach any copied vertex.  Hence
$S_{\mathrm{bi}}$ is $Z$-Grundy and has length $2r+2$.

\smallskip\noindent
\emph{Upper bound for the split graph.}
Take an arbitrary $Z$-Grundy sequence $T$ in $G_H^{\mathrm{sp}}$.
Once the first clique vertex has been selected, every vertex of
$B\cup\{c_0,c_1\}$ is closed-dominated.  Consequently all selected
copied vertices precede all selected clique vertices.  By
Lemma~\ref{lem:zcorebounds}, at most $r$ copied vertices are selected.
If the sequence contains $q\ge1$ clique vertices, all but the first
must have a fresh footprint in $A^*$, since the clique has already
been closed-dominated.  The same lemma gives $q-1\le r$.  Thus
$|T|\le2r+1$.  If $q=0$, the bound holds immediately.

\smallskip\noindent
\emph{Upper bound for the bipartite graph.}
Now we take any $Z$-Grundy sequence in $G_H^{\mathrm{bi}}$; selections
from the two bipartition classes may be mixed in any order.
By Lemma~\ref{lem:zcorebounds}, at most $r$ selected vertices are
copied vertices, and at most $r$ selected $B$-vertices have a fresh
footprint in $A^*$.
Every remaining selected vertex of $B$ must use $h$ as a fresh
footprint.  The same is true for $p_0$ or $p_1$.
Since a vertex can be used as a fresh footprint only once, the total
number of these exceptional selections is at most one.  Finally,
$h$ itself can be selected at most once.  Therefore
 $|T|\le r+r+1+1=2r+2.$
Combining each upper bound with its corresponding lower bound proves
\eqref{eq:zspidentity} and~\eqref{eq:zbiidentity}.
\end{proof}

\begin{theorem}\label{thm:znpcsplitbip}
\textsc{$Z$-Grundy Domination} is NP-complete on split graphs and on
bipartite graphs. 
\end{theorem}
\begin{proof}
Membership in NP follows because a proposed sequence can be checked in polynomial time.

For NP-hardness, take an instance $(H,k)$ of \textsc{Maximum Uniquely
Restricted Matching}.  We may delete isolated vertices and solve
edgeless instances directly.  Construct $G_H^{\mathrm{sp}}$ and
$G_H^{\mathrm{bi}}$ as above.  By Theorem~\ref{thm:zidentities},
\begin{align*}
 \nu_{\mathrm{ur}}(H)\ge k
 &\iff \gamma_{\mathrm{gr}}^Z(G_H^{\mathrm{sp}})\ge2k+1,\\
 \nu_{\mathrm{ur}}(H)\ge k
 &\iff \gamma_{\mathrm{gr}}^Z(G_H^{\mathrm{bi}})\ge2k+2.
\end{align*}
Both are polynomial-time reductions, and both target graphs
are connected.  NP-hardness follows from the source problem's
NP-completeness~\cite{golumbic2001uniquely}.
\end{proof}
\begin{corollary}
 \textsc{$Z$-Grundy Domination} is NP-complete  on word-representable graphs.  
\end{corollary}
\begin{proof}
    Hardness on word-representable graphs is immediate, as bipartite graphs are word-representable.
\end{proof}
The connection between $Z$-Grundy domination and standard zero
forcing is established by Bre\v{s}ar et al.~\cite{brevsar2017grundy}:
\begin{equation}\label{eq:zzfknown}
 Z(G)=|V(G)|-\gamma_{\mathrm{gr}}^Z(G).
\end{equation}
For the two constructions,
$|V(G_H^{\mathrm{sp}})|=2|A|+|B|+2$ and
$|V(G_H^{\mathrm{bi}})|=2|A|+|B|+3$.
Substitution into~\eqref{eq:zzfknown} gives the following expression:
\begin{equation}\label{eq:zforcingboth}
 Z(G_H^{\mathrm{sp}})=Z(G_H^{\mathrm{bi}})
 =2|A|+|B|+1-2\nu_{\mathrm{ur}}(H).
\end{equation}

\begin{corollary}\label{cor:zforcingboth}
\textsc{Minimum Standard Zero Forcing Set} is NP-complete on split
graphs, bipartite graphs as well as word-representable graphs. 
\end{corollary}
\begin{proof}
For either constructed graph, the identity~\eqref{eq:zforcingboth}
reduces the source instance $(H,k)$ to the zero forcing threshold
$2|A|+|B|+1-2k$.  Membership in NP follows by checking an initial
zero forcing set and a chronological sequence of forcing moves. 
\end{proof}

\section{NP-completeness on well-partitioned chordal graphs}\label{wpcs}

\begin{theorem}\label{thm:WPC}
Maximum Ordinary Grundy Domination is NP-complete on well-partitioned chordal graphs.
\end{theorem}
\begin{proof}
Bre\v{s}ar, Gologranc, Milani\v{c}, Rall, and Rizzi~\cite{brevsar2014dominating} reduce the NP-complete Grundy Covering Number problem on an arbitrary hypergraph $\mathcal H=(X,\mathcal E)$ to Grundy Domination on a chordal graph $G_{\mathcal H}$.  We recall the construction to verify that the resulting graph is also WPC.

The graph contains a central clique $K=X\cup\{x_0\}$, and, for every  hyperedge $F\in\mathcal E$, three vertices denoted $F$, $F^{(1)}$, and $F^{(0)}$.  Add the edges $Fx$ for $x\in F$, together with $FF^{(1)}$ and $F^{(1)}F^{(0)}$, and no other edges outside $K$.  The cited reduction proves the exact decision equivalence
\(
\rh(\mathcal H)\ge k
\quad\Longleftrightarrow\quad
\ggr(G_{\mathcal H})\ge 2|\mathcal E|+k+1.
\)
The source hypergraph may be assumed to have nonempty hyperedges, so every branch meets the central clique.

Now, for WPC graph, consider $K$ as a clique bag and for each hyperedge $F$, $\{F\},\{F^{(1)}\},\{F^{(0)}\}$ as singleton clique bags. Now, these bags are connected by the path
\(
K-\{F\}-\{F^{(1)}\}-\{F^{(0)}\}.
\)
All branches share the bag $K$, so the bag adjacency graph is a tree.  The interbag edges across the first edge are exactly $F\times\{F\}$, and the next two interbag edges join singleton bags.  Every interbag edge set is therefore complete bipartite between the prescribed boundary subsets.  This is a valid WPC partition tree and is computable in polynomial time.

\end{proof}

\section{Conclusion}
\label{final}

In this paper, we studied Grundy domination and its
variants on co-bipartite, split, bipartite, and
well-partitioned chordal (WPC) graphs. We first
developed polynomial-time algorithms for finding
maximum Grundy covering sequences in interval and
circular-arc hypergraphs. Using these algorithms,
we obtained polynomial-time algorithms for ordinary,
$Z$-, total, and $L$-Grundy domination on
circular-convex co-bipartite graphs. We also
established an equality between the ordinary and
$Z$-Grundy domination numbers on connected
co-bipartite graphs.

On the complexity side, we proved that the
$Z$-Grundy Domination problem is NP-complete on
both split and bipartite graphs through reductions
from the Maximum Uniquely Restricted Matching
problem. These reductions also yield NP-completeness
results for the standard zero forcing problem.
Furthermore, we established NP-completeness of
ordinary Grundy domination on WPC graphs,
even when a valid partition tree is supplied.


\end{document}